\documentclass[10pt]{article}

\usepackage{epsf,epsfig,amsfonts,graphicx,a4wide} 
\usepackage{amssymb}
\usepackage{amsthm}
\usepackage{enumerate}
\usepackage{url}
\usepackage{epsf,amsfonts,graphics}
\usepackage{graphicx} 
\usepackage{color}
\usepackage{a4wide,epsfig}
\newtheorem{theorem}{Theorem}

\newtheorem{remark}{Remark}

\newcommand{\weg}[1]{}

\newcommand{\Id}{\textrm{Id}}

\title{ Zermelo deformation of Finsler metrics  by   Killing vector fields  } 
\author{   Patrick  Foulon and  
Vladimir  S. Matveev  }

\date{} 

\begin{document}

\maketitle
 
 \begin{abstract} 
We show how geodesics, Jacobi vector fields and flag curvature of a Finsler metric behave under 
 Zermelo deformation with respect to a Killing vector field. We also show that Zermelo deformation with respect to a Killing vector field of a locally symmetric Finsler metric is also locally symmetric. 
\end{abstract} 

   \section{Introduction}  
\begin{figure}[ht]
	\centering
		\includegraphics[width=.5\textwidth]{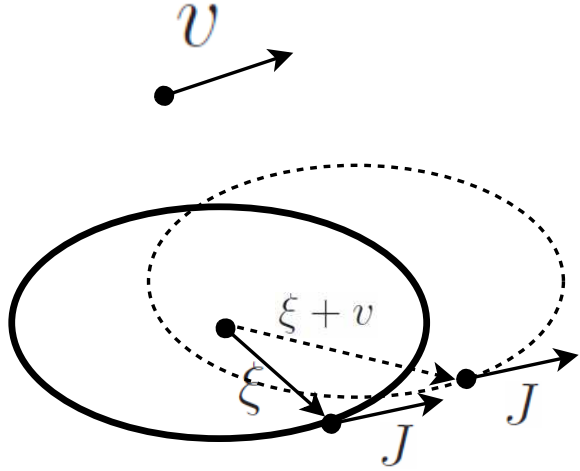}
	\caption{The unit ball of $\tilde F$ (punctured line)  is the $v$-translation of that of $F$ (bold line). If a vector $J$ is tangent to the  unit ball of $ F$ at $\xi$, it is tangent to  the  unit ball of $ \tilde F$ at $\xi + v$} \label{fig:randers} 
\end{figure}

 Let $F$ be a Finsler metric on $M^n$ and  $v$ be  a vector field such that $F(x,-v(x))<1$ for any $x\in M^n $. 
 We will denote by $\tilde F$ the Zermelo deformation of $F$ by   $v$. That is, for each point $x\in M$,  the unit $\tilde F$-ball  $\tilde B_x:= \{\xi\in T_xM^n \mid \tilde F(x,\xi)<1\}$  is  the translation
in  $T_xM^n$  along the vector $v(x)$ of the unit   $F$-ball 
 $B_x:= \{\xi\in T_xM^n \mid  F(x,\xi)<1\}$ (see Fig. \ref{fig:randers}).  

Equivalently, this  can be reformulated as

   \begin{equation} \label{tildeF} 
  \tilde F(x,\xi) =  F(x, \xi - \tilde F(x,\xi) v(x)).\end{equation}  Indeed,  the equation (\ref{tildeF}) is positively homogeneous, and for any $\xi$ such that $\tilde F(x,\xi)=1$  we have $F(x,\xi -v(x))= 1$.

  The  first  result   of this note is a description of   how geodesics, Jacobi vector fields  and flag curvatures  of 
$F$ and of $\tilde F$ are related, if the vector field $v$ is a Killing vector field for $F$, that is, if the flow of $v$ preserves $F$.

\begin{theorem} \label{flag}  Let $F$ be a Finsler metric on $M^n$  
 admitting a Killing vector field $v$ such that $F(x,-v(x))<1$ for all $x\in M^n$. We denote by  $\Psi_t$ the  flow of $v$  and by $\tilde F$ the $v$-Zermelo deformation of $F$.

 Then, for any $F$-arc-length-parameterized geodesic $\gamma$ of $F$, the curve $t\mapsto \Psi_t (\gamma(t))$ which we denote by $\tilde \gamma(t)$ is a $\tilde F$-arc-length-parameterized geodesic   of $\tilde F$.

 Moreover,  for any   Jacobi vector field $J(t)$ along $\gamma$ such that it is orthogonal to 
$\dot \gamma(t)$ in the metric  $g_{(\gamma(t), \dot \gamma(t))}:= \frac{1}{2} d_{\xi}^2F^2_{(\gamma(t), \dot\gamma(t))}$,   
 the  pushforward  $\tilde J(t)= \Psi_{t*}(J(t))$ is  a Jacobi vector field for  $\tilde \gamma(t)$ and   is  
  orthogonal to $\dot{\tilde \gamma}(t)$ in  $\tilde g_{(\tilde \gamma(t), \dot{\tilde \gamma}(t))}:=  \frac{1}{2} d_{\xi}^2\tilde F^2_{(\tilde\gamma(t), \dot{\tilde \gamma}(t))}.$

  Moreover, flag curvatures $K$ and $\tilde  K$ of $F$ and $\tilde F$ 
  are related by the following formula: for any $x\in M$ and any ``flag'' $(\xi, \eta )$  with 
flagpole $\xi\in  T_xM $ and transverse  edge $\eta\in  T_xM $, we have $K(x, \xi, \eta)=  \tilde K(x, \xi + v, \eta)$ provided that $\xi + v$ and $ \eta$ are linearly independent.  
\end{theorem}

  We do not pretend that the whole result  is new and  rather suggest  that its certain parts are  known. The first statement of Theorem \ref{flag} appears in \cite{Katok}. We recall  the arguments of A. Katok in Remark \ref{kat}. The third  statement was announced in \cite{Foulon} and follows from the recent paper \cite{HuangMo}.   
  Special cases when the metric $F$ is Riemannian were
  studied in details in  e.g.  \cite{BRS, Shen}. Though we  did not find the second statement of Theorem \ref{flag}, the one about the Jacobi vector fields, in the literature, we think it is  known in folklore.

 Unfortunately in all these references, the proof is by direct calculations, which are sometimes quite tricky and sometimes require a lot of preliminary work. One of the  goals of this note is to show the geometry lying below Theorem \ref{flag}
  and to demonstrate that certain parts of Theorem \ref{flag}
  at least are  almost trivial.

 Our second result shows that Zermelo deformation with respect to a Killing vector field preserves the property of a Finsler metric to be a locally symmetric space.  We will call  a Finsler metric {\it locally symmetric}, if for any geodesic $\gamma$ the covariant 
 derivative of the Riemann curvature (=Jacobi operator)  vanished:
 \begin{equation} \label{cov} 
 D_{\dot \gamma} R_{\dot \gamma}= 0. 
 \end{equation}
 Here $D_{\dot \gamma}$  stays for the covariant derivative along the geodesic: $D_{\dot \gamma} = \nabla^{\dot \gamma}_{\dot \gamma} $. Both most popular Finslerian  connections, Berwald and Chern-Rund connections,  can be used as the Finslerian  connection in the last formula, see e.g. \cite[\S 7.3]{shen2}, whose notation we partially follow.

\begin{remark} \label{kat} 
In Riemannian geometry there exist
 two equivalent definitions of locally  symmetric spaces: according to the ``metric''  definition, a space is locally symmetric  if for   any point there exists a local  
 isometry such that this point is a fixed point and the differential of the isometry at this point  is minus the identity. By the  other  ``curvature'' definition,   a space is locally symmetric if  the covariant derivative of the curvature tensor is zero. The equivalence of these two definitions is a classical result of E. Cartan. We see that our definition above is the generalization to the Finsler metrics of the ``curvature'' definition; it  was first suggested in \cite{foulon}. 
 
 In the  Finsler setup, both ``metric'' and ``curvature''  definitions are used in the literature, but they  are not  equivalent anymore: the symmetric spaces with respect to the ``metric'' definition are symmetric spaces with respect to the ``curvature'' definition, but not vice versa. 
 
 In fact, the ``metric'' definition is much more restrictive; in particular locally symmetric metrics  in the ``metric'' definition are  automatically Berwaldian \cite[Theorem 9.2]{MT2012} and are clearly reversible. On the other hand, all metrics of constant flag curvature,  in particular all Hilbert metrics in a  strictly convex domain,  are locally  symmetric  in the ``curvature''  definition,  and are symmetric with respect to the ``metric''  definition if and only if the domain is an ellipsoid. 

In view of this, the name ``locally symmetric'' is slightly misleading, since  locally symmetric manifolds may  have no  (local) isometries.  We will still use this terminology because it was used in the literature before.   
 \end{remark}

 \begin{theorem} \label{thm:2}
 Suppose that $F$  is a locally symmetric Finsler metric  and $v$ a Killing vector field satisfying $F(x,-v(x))<1$ for all  $x$. Then, the $v$-Zermelo deformation of $F$ is also locally symmetric. 
 \end{theorem}

 All Finsler metrics in our paper are assumed to be   smooth and strictly convex but may be  irreversible.

 \section{Proofs.}

\subsection{Proof of Theorem \ref{flag}.}  Let $\gamma(t)$ be  an arc-length-parameterized  
 $F$-geodesic, we need to prove that  the curve  $t\mapsto \Psi_{t}(\gamma(t))$  is   an arc-length-parameterized  
  $\tilde F$-geodesic.  
   In order to do it,  observe that 
for any $F$-arc-length-parameterized  curve  $x(t)$ the  $t$-derivative of $\Psi_t(x(t))$ is given by $\Psi_{t*}(\dot x(t)) + v(\Psi_t(x(t)))$. Since  the flow 
$\Psi_{t}$ preserves $F$ and $v$, it preserves $\tilde F$ and therefore  
$$\tilde F\left(\Psi_t(x(t)), \Psi_{t*}(\dot x(t)) + v(\Psi_t(x(t)))\right)= \tilde F\left(x(t), \dot x(t)  + v(x(t))\right)=  F\left(x(t), \dot x(t) \right).$$ The last equality in the formula above  is true because, for any $\xi $ such that $F(x, \xi)= 1$, we have  
$ \tilde F\left(x  ,\xi  +v(x)  \right) = F(x,\xi)$  by the definition of the Zermelo deformation.     Thus,  if the curve  $x(t)$ is $F$-arc-length-parameterized, then  the curve  $\Psi_t(x(t))$ is $\tilde F$-arc-length-parameterized. 
   
 This also 
  implies that the 
 integrals $\int F(x(t), \dot x(t)) dt$ and $\int \tilde F\left( \Psi_t(x(t)), \dot{\left(\Psi_t(x(t))\right)} \right)dt $ coincide for all $F$-arc-length-parameterized curves $x(t)$. Since geodesics are locally  the shortest arc-length parameterized curves connecting two points, 
     for each arc-length parameterized $F$-geodesic $\gamma$  the curve $t\mapsto \Psi_t(\gamma(t))$ is a $\tilde F$-arc-length-parameterized geodesic  as we claimed.

\begin{remark}  Alternative geometric proof of the statement that  for each arc-length parameterized $F$-geodesic $\gamma$  the curve $t\mapsto \Psi_t(\gamma(t))$ is 
a $\tilde F$-arc-length-parameterized geodesic is essentially due to \cite{Katok}: consider the Legendre-transformation $T:T^*M^n \to TM$ 
 corresponding to the function $\frac{1}{2} F^2$ and denote by $F^*$ the pullback of $F$ to $T^*M^n$, \ $F^*:= F \circ T$. Next, 
  view the vector field 
 $v$   as a function on $T^*M$ by the obvious rule $\eta\in T^*_xM^n \mapsto \eta(v(x))$. It is known that the Hamiltonian flow corresponding to the  function $v$ is the natural lift of the flow of the vector field $v$ to $T^*M$. Since $v$ is assumed to be a Killing vector field,  the Hamiltonian flows of $F^*$ and of $v$ commute. Next, consider the function $\tilde F^*:= F^* +v$. If $v$ satisfies 
$F(x,-v(x))<1$, then  the restriction of $\tilde F^*$ to $T_xM^n $  is convex, consider the  Legendre-transformation $\tilde T:TM^n \to T^*M$ corresponding to the function $\frac{1}{2} (\tilde F^*)^2$ and the pullback of $\tilde F^*$ to $TM^n$, it is a Finsler metric which we denote by $\tilde F$. It is a standard fact in convex geometry that the Finsler metric $\tilde F$ is the $v$-Zermelo-deformation of $F$.  Since the Hamiltonian flows of $F^*$  and of $v$, which we denote by $\psi_t$ and $ d^* \Psi_t$ , 
 commute, the Hamiltonian flow  of $\tilde F^* $ is simply given by  \begin{equation} \label{com} 
 \tilde \psi_t =   d^* \Psi_t \circ\psi_t.\end{equation}
  Then, for any point  $(x, \xi)\in TM$ with $F(x,\xi)=1$,  the projections of the  orbits of 
  $\tilde \psi_t$ and of $\psi_t$   starting  at this  point  are arc-length parameterized  geodesics $\gamma$ of $F$ and $\tilde \gamma$ of $\tilde F$. By (\ref{com}) we have 
    $\tilde \gamma(t)=  \Psi_t\circ \gamma(t)$ as we claimed. 
 
\end{remark}

Let us now  prove the second statement of Theorem \ref{flag}. Consider a  Jacobi vector field $J(t)$   which are orthogonal to $\gamma$.   
We need to show that  the pushforward  $\tilde J(t)= \Psi_{t*}(J(t))$ is  a Jacobi vector field for the $\tilde F$-geodesic  $\tilde \gamma(t):= \Psi_{t}(\gamma(t))$. By the definition of Jacobi vector field there exists 
 a family $\gamma_s(t)$ of geodesics  with $\gamma_0= \gamma$ such that $J(t)= \frac{d}{ds}_{|s=0} \gamma_s(t)$, since $J(t)$   is orthogonal to $\gamma$  we may assume that all geodesics $\gamma_s(t)$ are arc-length parameterized.  As we explained above,  $\Psi_{t}(\gamma_s(t))$ is a family   of $\tilde F$-geodesics; taking the derivative by $s$ at $s=0$ proves what we want.

Let us  now show  that $\tilde J$ is orthogonal to $\dot{\tilde\gamma}$. 
First observe that  the condition that $J(t)$ is orthogonal to $\dot \gamma(t)$  is  equivalent to the condition that 
  $J(F):= \sum_r J^r \frac{\partial F}{\partial \xi_r} $  vanishes at  $\dot \gamma(t)$ for each $t$. Indeed, consider the one-form  
	$U\in T_{\gamma(t)}M^n \mapsto g_{\gamma(t), \dot \gamma(t)}(\dot \gamma(t)_,U  )$. 
	Because of the (positive) homogeneity of the function  $F$ we have that at a point $\dot \gamma(t)\in T_{\gamma(t)}M^n$   
 \begin{equation} \label{3}  g_{(\gamma(t), \dot \gamma(t))}(\dot \gamma(t)_,U  ) =    U^r  \frac{\partial F}{\partial \xi_r}.   \end{equation} 
   
Next, take   Equation (\ref{tildeF}) and  calculate the differential of  the restriction of 
  $\tilde F$ to the tangent space: its  components are given by 
   \begin{equation} \label{dtildeF} 
  \frac{\partial \tilde F }{\partial \xi_i}= \frac{1}{1 + v(F)} \frac{\partial F }{\partial \xi_i}. 
\end{equation}
 In this formula, the derivatives of the function $ F$ are  taken at  $\xi\in T_xS^2$, and the derivatives of the function $\tilde F$ are  taken at  $\xi- \tilde F(x,\xi) v$. By $v(F)$ we denoted the function $\sum_{r} \frac{\partial F}{\partial \xi_r} v^r$. 

 In view of (\ref{dtildeF}), 
$ J( \tilde F):= \sum_r  J^r \frac{\partial \tilde F}{\partial \xi_r} $  vanishes at  $\dot{ \gamma}(t) + v(\gamma(t))$,  so 
$ J(t) $ is orthogonal to $\dot {\gamma}(t) + v(\gamma(t))$ (the orthogonality is understood in the sense of $\tilde g_{(\gamma(t), \dot \gamma(t) + v(\gamma(t)))}$). Then, $\tilde J(t)= \Psi_{t*}(J(t))$ is $\tilde g_{(\tilde \gamma(t), \dot {\tilde \gamma}(t))}  $ 
orthogonal to $\dot{\tilde \gamma}(t).$

\begin{remark} 
Geometrically, the just proved 
statement that  $\tilde J$ is orthogonal to $\dot{\tilde \gamma}$,  after the identification of $T_{\gamma(t)}M^n $ and $T_{\Psi_t(\gamma(t))}M^n$ by the differential of the diffeomorphism  $\Psi_t$,  
 corresponds to the following simple observation:   if $J$ is tangent to the unit $F$-sphere at the point $\xi=\dot \gamma$, then it is also tangent to the unit $\tilde F$-sphere at the point $\xi + v$, see Fig. \ref{fig:randers}.   
\end{remark}

 Let us now prove the third statement of Theorem \ref{flag}, we need to show that $K(x, \xi, \eta)= \tilde K(x, \xi+ v, \eta)$. We consider a $F$-geodesic 
 $\gamma(t)$ with $\gamma(0)= 0$ and $\dot \gamma(0)= \xi$ and the corresponding $\tilde F$-geodesic $\tilde\gamma(t)= \Psi_t(\gamma(t)); $ since $\Psi_0= \Id$,  we have  $\dot{\tilde \gamma}(0)= \dot \gamma(0)+ v = \xi + v$.   
  
 Observe that by
  combining \cite[Eqn. 6.16 on page 117]{nankai}  and  \cite[Eqn. 6.3 on page 108]{nankai} we obtain 
\begin{equation} \label{skoro} \begin{array}{lcl}
\frac{1}{2}\frac{d^2 }{dt^2} g(J(t), J(t))&=& g(D_{\dot \gamma(t)} D_{\dot \gamma(t)} J(t), J(t)) + g(D_{\dot \gamma(t)}J(t),  D_{\dot \gamma(t)} J(t)) \\  & =& -  K (\gamma(t), \dot \gamma(t), J(t)) \left(g(\dot\gamma(t), \dot\gamma(t)) g(J(t), J(t))  -  g(\dot\gamma(t), J(t))^2  )  + g(D_{\dot \gamma}J,  D_{\dot \gamma} J\right).  \end{array}\end{equation}

We will assume that $J(0):= \eta$ is $g$-orthogonal to $\dot \gamma(0)$. As explained above this implies that $\dot {\tilde \gamma}(0)$ is $\tilde g$-orthogonal to $\tilde J(0)$. 
Then, by (\ref{skoro}),   the minimum of $\frac{d^2 }{dt^2} g(J(t), J(t))_{|t=0}$ 
 taken over  all Jacobi vector fields $J$  along $\gamma$ 
  which  are equal to $\eta$ at $t=0$,    is equal to $-K(x, \xi, \eta) g(\xi, \xi) g(\eta, \eta) $. 
   An analogous statement  is clearly true also for $\tilde \gamma$, $\tilde g$ and 
  $ \tilde J$:  namely, 
  the minimum of $\frac{d^2 }{dt^2} \tilde g(\tilde J(t), \tilde J(t))_{|t=0}$ 
 taken over  all Jacobi vector fields $\tilde J$  along $\tilde \gamma$ 
  which  are equal to $\eta$ at $t=0$,   is equal to $-\tilde K(x, \xi + v  , \eta)\tilde g(\xi + v, \xi + v) \tilde g(\eta, \eta) $. Here we used the relation  
   $\tilde J(0)= \Psi_{t*}(J(t))_{|t=0} = J(0)$. Finally,    in order to show that  $K(x, \xi, \eta) = \tilde K(x, \xi+ v, \eta),  $
 it is sufficient to show that 
  the function $t\mapsto g(J(t), J(t))$  is proportional, with a constant coefficient, to  the function 
 $t\mapsto  \tilde g(\tilde J(t), \tilde J(t)). $

 In order to prove this,   let us 
 first compare   $g_{(\gamma(t), \dot \gamma(t))}= \frac{1}{2} d_{\xi}^2F^2_{(\gamma(t), \dot\gamma(t))} $ and 
  $$\tilde g_{(\tilde \gamma(t), \dot {\tilde \gamma}(t)  )}= \frac{1}{2} d_{\xi}^2\tilde F^2_{(\Psi_t \circ \gamma(t), \Psi_{t*} (\dot \gamma(t)) + v(\Psi_t(\gamma(t))))}. $$
		It is convenient to work in coordinates $(x_1,...,x_n)$ such that  the entries of $v$
	are constants,  in these coordinates for each $t$  the differential of the diffeomorphism $\Psi_t$ is given  by the identity matrix, so in  these  coordinates  $J(t) =  \tilde J(t)$  and $  
  \dot {\tilde \gamma}(t)=\dot \gamma(t)+ v(\gamma(t)).$

 Differentiating (\ref{dtildeF}), we get  the second derivatives of $\tilde F$. They are given by 
\begin{equation} \label{tildegg} 
  \frac{\partial^2 \tilde F }{\partial \xi_i \partial \xi_j} = \frac{1}{1 + v(F)} \frac{\partial^2 F }{\partial \xi_i \partial \xi_j} -
 \frac{1}{(1 + v(F))^2} \sum_r\left( \frac{\partial^2 F }{\partial \xi_i \partial \xi_r} v^r  \frac{\partial F }{\partial \xi_j} + 
   \frac{\partial^2 F }{\partial \xi_j \partial \xi_r} v^r  \frac{\partial F }{\partial \xi_i}\right). 
\end{equation} 
 Again, all  derivatives of the function $F$ are  taken at  $\xi$, and  of the function $\tilde F$ are  taken at  $\xi- \tilde F(x,\xi) v(x)$. Note that one term in the brackets in (\ref{tildegg})
 appears   because we differentiate $\frac{1}{1 + v(F)}$, and the other  appears  because the derivatives  of 
  $ \frac{\partial F }{\partial \xi_j}$ are taken at  $\xi- \tilde F(x,\xi) v(x)$. When we differentiate it, we also need to take into account the additional term $- \tilde F(x,\xi) v(x)$.

 Now, in view of the formula $\frac{1}{2}d^2(\tilde F^2)= \tilde F d^2\tilde F +   d\tilde F \otimes d\tilde F$ we   obtain  from  (\ref{tildegg}) the formula for $\tilde g_{ij}$:
  \begin{equation} \label{tildeg} 
 \tilde g_{ij} = \frac{\tilde F}{1+ v(F)}  g_{ij}-  \frac{\tilde F}{(1 + v(F))^2}\sum_r\left( \frac{\partial^2 F }{\partial \xi_i \partial \xi_r}  v^r  \frac{\partial F }{\partial \xi_j} + 
 \frac{\partial^2 F }{\partial \xi_j \partial \xi_r}  v^r  \frac{\partial F }{\partial \xi_i}\right) + \frac{1}{(1 + v(F))^2}\frac{\partial   F }{\partial \xi_i}\frac{\partial   F }{\partial \xi_j}
.  \end{equation} 
 Let us now compare the length of $J$ in $g(x,\xi)$ with that of  in $\tilde g(x, \xi +v)$. We multiply (\ref{tildeg}) by $J^iJ^j$ and sum  with respect to $i$ and $j$. Since by assumptions   $J(F)= \sum_r J^r \frac{\partial  F}{\partial \xi_r} $ vanishes at $\xi$, all terms in the sum but the first vanish. We thus obtain that the length of $J$ in $\tilde g$ is proportional to that of in $g$ with the coefficient which is 
the square root of  $\frac{\tilde F}{1+ v(F)}$. 

But along the geodesic both $\tilde F$ and $v(F)$ are constant. 
Indeed, $v(F)$ is the ``Noether'' integral corresponding to the Killing vector field. 
 Theorem \ref{flag} is proved.

\subsection{Proof of Theorem \ref{thm:2}.} 
 First, observe that a Finsler metric is locally symmetric if and only if for  any geodesic 
$ \gamma$ and  any Jacobi  vector field $  J$ along $\gamma$   
the vector field 
$ D_{\dot  \gamma}  J $   is also a Jacobi vector field. Indeed, the equation for Jacobi vector fields is 
\begin{equation} \label{Jac} 
D_{\dot \gamma} D_{\dot \gamma} J + R_{\dot \gamma}(J)= 0. 
\end{equation}   
$D_{\dot \gamma}$-differentiating this equation, we obtain 
$$
D_{\dot \gamma} (D_{\dot \gamma} D_{\dot \gamma} J + R_{\dot \gamma}(J)) =  D_{\dot \gamma} D_{\dot \gamma} (D_{\dot \gamma} J) + R_{\dot \gamma}(D_{\dot \gamma} J) + \left(D_{\dot \gamma}R_{\dot \gamma}\right)(J)=0.
$$
If $D_{\dot \gamma} J$ is a Jacobi vector field, 
$D_{\dot \gamma} D_{\dot \gamma} (D_{\dot \gamma} J) + R_{\dot \gamma}(D_{\dot \gamma} J)$ vanishes  so the  equation above implies $\left(D_{\dot \gamma}R_{\dot \gamma}\right)(J)=0$, and since it is fulfilled for all Jacobi vector fields we have 
$D_{\dot \gamma}R_{\dot \gamma}=0$ as we claimed.

Thus, we assume that for any geodesic and for 
any Jacobi vector field for $F$ its  $D_{\dot \gamma}$ derivative is also a Jacobi vector field, and our goal is to show the same for $\tilde F$. Clearly, it is sufficient to show this only for Jacobi vector fields which are $g$-orthogonal to $\dot { \gamma}$.  Note that for such Jacobi vector fields  $ D_{\dot { \gamma}}  J$ 
 is also  orthogonal to $\dot \gamma$, since both  $g_{(\gamma, \dot \gamma)}$ and $\dot\gamma$ are 
			$ D_{\dot \gamma}$-parallel.

Take a (arc-length-parameterized) $F$-geodesic $\gamma$ and  a point $P= \gamma(0)$
 on it.  Consider the geodesic 
polar coordinated around this point, let us recall what they are and their properties which we use in the proof. 

 Consider   the  (local) diffeomorphism of $T_PM^n \setminus \{0\}$ to $M^n $ which sends    $\xi \in T_PM^n\setminus \{0\} $ to  $\exp(\xi):= \gamma_{\xi}(1)$, where $\gamma_\xi$ is the geodesic starting from $P$ with the velocity vector $\xi$. As the local coordinate systems on $T_PM^n\setminus \{0\} $ we take the following one: we choose a local coordinate system $x_1,...,x_{n-1}$ 
on the unit $F$-sphere $\{\xi \in T_PM^n \mid F(P, \xi)= 1\}$  and set 
the tuple  $\left(F(P, \xi), x_1\left(\frac{1}{F}\xi\right),...,  x_{n-1}\left(\frac{1}{F}\xi\right)\right)$ to be the coordinates of $\xi$.  Combining it with the diffeomorphism $\exp$, we obtain a local coordinate system on $M^n$. By construction,  in this coordinate system each  arc-length parameterized geodesics  starting at $P$, in particular the geodesic $\gamma$,
 is  a curve  of the form  $(t, \textrm{const}_1,..., \textrm{const}_{n-1})$.

 Next, consider the following local Riemannian metric $\hat g$  in a punctured neighborhood of $P$: for a point $\sigma(t)$ of this neighborhood such that $\sigma$ is a geodesic passing through $P$ we set $\hat g:= g_{(\sigma(t),\dot \sigma(t))}$.   
 
It is known that in the polar coordinates the metric $\hat g$ is block-diagonal with one  $1\times 1$ block  which is simply the identity
and one $(n-1)\times (n-1)$-block which we denote by $G$:
$$
\hat g = \left(\begin{array}{cc} 1 & \\ & G\end{array}\right).$$  
It is known, see e.g. \cite[Lemma 7.1.4]{shen2}, that 
the geodesics passing through $P$ are also geodesics of $\hat g$, and that for each such geodesic the operator $\hat 
\nabla_{\dot \gamma}$, where $\hat \nabla$ is the Levi-Civita connection of $\hat g$, coincides with $ D_{\dot \gamma}$.

Next, consider analogous objects for  the  metric $\tilde F$. As the local  coordinate system on the unit $\tilde F$-sphere we take the following: as the coordinate tuple 
 of   $\tilde \xi$ with $\tilde F(P, \tilde \xi)=1$ we take   the coordinate tuple $(x_1( \xi),..., x_{n-1}(\xi))$, where $\xi:=\tilde \xi -v(P)$. 
   (Recall that  the  $v$-parallel-transport  sends  the unit $F$-sphere  to the unit $\tilde F$-sphere.) 
   
   By Theorem \ref{flag},  in these coordinate systems    each Jacobi vector field 
   $J= (J_0(t),..., J_{n-1}(t))$
    along $\gamma$  which is orthogonal to $\gamma$   is also a Jacobi vector field along $\tilde \gamma$,  which is the $\tilde F$-geodesic such that  $\tilde \gamma(0)= P$ and  $\dot{\tilde  \gamma}(0) = \dot \gamma(0) +  v(P)$, and is orthogonal to $\tilde \gamma$.   
   By (\ref{tildeg}), the corresponding block $\tilde G$ is given by $\tilde G =\frac{1}{1 + v(F)} G$. Since the function $v(F)$ is constant along geodesics, the coefficients $\Gamma_{ij}^k$ 
   of the Levi-Civita connection $\hat \nabla$ of $\hat  g$  such that  $i = 0$ or   $j = 0$  coincide with that of for the analog for  $\tilde F$. A direct way to see the last claim is to use the formula $\Gamma_{ij}^k= \frac{1}{2} g^{ks} \left(\frac{\partial g_{is}}{\partial x_j} + \frac{\partial g_{js}}{\partial x_i} - \frac{\partial g_{ij}}{\partial x_s}  \right) $, where of course all indices run from $0$ to $n-1$ and the summation convention  is assumed.

   Then, in our chosen coordinate system, 
   the formula for the covariant derivative in $\hat \nabla$ along $\gamma$ 
   for vector fields which are orthogonal to $\gamma$ simply coincides with that of the formula for the corresponding objects for 
   $\tilde F$. 
    Then,  for any $\tilde F$-Jacobi vector field $\tilde J$  orthogonal to $\dot {\tilde \gamma}$ 
    we have that 
      $\tilde D_{\dot {\tilde \gamma}} \tilde J$- is again a Jacobi vector field.     Theorem \ref{thm:2} is proved.

\vspace{1ex}
{\bf Acknowledgments.}
The authors thank  Sergei Ivanov for useful comments. V.M.  was partially supported  by  the University of  Jena and by  the DFG grant MA 2565/4.

Patrick Foulon:\ Centre International de Rencontres Math\'ematiques-CIRM, 163 avenue de Luminy,
Case 916, F-13288 Marseille - Cedex 9, France. \\ {\bf Email:  \url{foulon@cirm-math.fr}}

\vspace{1ex} 

 Vladimir S. Matveev: {Institut f\"ur Mathematik,
Fakult\"at f\"ur Mathematik und Informatik,\newline
Friedrich-Schiller-Universit\"at Jena,
07737 Jena, Germany. \\ {\bf Email: \url{vladimir.matveev@uni-jena.de}  }

\end{document}